\documentclass[11pt]{amsart}
\usepackage{latexsym,amssymb,amsmath}
\newtheorem{theorem}{Theorem}[section]

\newtheorem{corollary}[theorem]{Corollary}

\theoremstyle{definition}

\theoremstyle{remark}

\numberwithin{equation}{section}

\def\CP{\mathbb{CP}}
\def\C{\mathbb{C}}

\def\gS{\mathfrak{S}} 
\def\l{\lambda}
\def\w{\omega}
\def\a{\alpha}

\def\d{\delta}

\def\r{\rho}

\def\frp#1{\frac{\partial}{\partial{#1}}}

\def\xx#1{x^{#1}}

\def\ooo#1#2{\omega^{#1}_{#2}}
\def\oo#1{\omega^{#1}_0}

\def\BC{\mathbb C}\def\BF{\mathbb F}
\def\BR{\mathbb R}
\def\BG{\mathbb G}
\def\BP{\mathbb P}
\def\pp#1{\mathbb P^{#1}}

\def\pp#1{{\mathbb P}^{#1}}
\def\tdim{\rm dim}
\def\hd{,...,}

\def\upperp{{}^\perp}

\def\inv{{}^{-1}}

\def\cC{{\mathcal C}}
\def\cF{{\mathcal F}}

\def\cO{{\mathcal O}}

\def\11{\mathbf 1}
\def\PP{\mathbb P}

\def\l{\lambda}
\def\a{\alpha}

\def\o{\omega}

\def\d{\delta}

\def\ot{{\mathord{\,\otimes }\,}}
\def\op{{\mathord{\,\oplus }\,}}

\def\ra{{\mathord{\;\rightarrow\;}}}

\def\tdim{\text{dim}\,}
\def\tcodim{\text{codim}\,}

\def\tmod{\text{ mod }}

\def\tmax{\text{ max }}

\begin{document}
\title{Fubini's theorem in codimension two}
\author{J.M. Landsberg 
 \and Colleen Robles}
\date{September, 2005}
\begin{abstract} We classify codimension two analytic submanifolds
of projective space $X^n\subset \BC\pp{n+2}$
having the property that any line through a general
point $x$ having contact to order two with $X$ at $x$ automatically has
contact to order three.  We give   applications
to the study of the Debarre-de Jong conjecture and
  of
  varieties  whose Fano variety of
lines has dimension $2n-4$.
\end{abstract}
\thanks{Landsberg supported by NSF grant DMS-0305829}
\email{jml@math.tamu.edu, robles@math.rochester.edu}
\maketitle


\section{Introduction}

\subsection{Statement of the main result}
Let $V$ be a complex vector space, and 
$X\subset \BP V$ be a complex submanifold or algebraic variety
and let $x\in X$ be a smooth point. Define $\cC_{k,x}\subset \BP T_xX$
to be the set of tangent directions at $x$ for which there
exists a line $l\simeq \pp 1$ in $\BP V$   having contact
to order $k$ with $X$ at $x$, or, in 
the language of algebraic geometry, $mult(l\cap X)_x\geq k+1$.
Let $\cC_x=\cC_{\infty, x}\subset \BP T_xX$ denote the tangent directions
to lines on $X$ through $x$.

One way to state the classical Fubini theorem \cite{fubini}
is as follows:

\begin{theorem}[Fubini] 
\label{thm:fubini}
Let  $X^n\subset \BC\pp{n+1}$ be a complex analytic hypersurface  with $n>1$ 
and at least a two dimensional Gauss image. Let $x\in X$ be a general point. 
If $$\cC_{2,x}=\cC_{3,x} \ \ {\rm (Fubini \ hypothesis)}$$
then $X$ is (an open subset of) a quadric hypersurface.
\end{theorem}

We stated the redundant hypotheses $n>1$ for emphasis. When $n=1$
the Fubini hypothesis is vacuous.
If  $X$ is a hypersurface whose Gauss image has dimension one, then 
$X$ is locally ruled by $\pp{n-1}$'s \cite{sev,BS48}. (I.e. if $X$ 
is variety, it is a scroll of $\pp{n-1}$'s.)  So, all hypersurfaces 
satisfying the Fubini hypothesis are classified.

\medskip

In this paper we present a generalization of Fubini's theorem
to codimension two. There are several formulations
of the Fubini hypothesis, all of which are equivalent for
hypersurfaces but do not all coincide already in codimension two.
Thus our first task is to come up with proper hypotheses.
Let $X\subset \BP V$ be a variety or analytic submanifold
and let $x\in X$ be a general point. There is a well
defined sequence of ideals defined on the tangent space
$T_xX$ given by the {\it relative differential invariants}
$F_k\in S^kT^*_xX\ot N_xX$, where $F_k$ is 
(an equivalence class of) vector spaces of homogeneous polynomials of degree
$k$ on $T_xX$ parametrized by the conormal space $N^*_xX$.
A coordinate definition of these invariants is as follows:
Take adapted local coordinates $( w^{\a},z^{\mu})$, $1\leq \a\leq n$,
$n+1\leq \mu\leq \tdim \BP V$,
on $\BP V$ such that $[x]=(0,0)$ and
$T_{[x]}X$ is spanned by the first $n$ coordinates ($1\leq \a\leq n$). Then locally $X$ is given by equations
\begin{equation}\label{grapheqn}
z^{\mu}=f^{\mu}(w^{\a})
\end{equation} and, at $(0,0)$,
$$
F_{k}(\frp {w^{i_1}}\hd \frp{w^{i_k}}) =\sum_{\mu}\frac {\partial^kf^{\mu}}
{\partial w^{i_1} \hd \partial w^{i_k}} \frp {z^{\mu}}.
$$
The invariant $II=F_2$ is called the {\it projective second fundamental
form} and for it there is no equivalence to mod out by. For the other
invariants, different choices, e.g.,  of a complement to $T_xX$ in
$T_x\BP V$, will yield different systems of polynomials, but
the new higher degree polynomials will be the old plus polynomials
in the ideal generated by the lower degree forms
(see \cite{IvL}, \S 3.5). Letting
$|F_k|=F_k(N^*_xX)\subseteq S^kT^*_xX$,
the ideals in $Sym(T^*_xX)$ generated by $\{|F_2|\hd |F_k|\}$
are well defined. 

The set $\cC_{k,x}$ is   the zero set of
$\{|F_2|\hd |F_k|\}$. Because points can and do occur with multiplicities,
it will be more precise to work with
the ideals $I_{\cC_{k,x}}$ which we define to be the ideals
generated by  $\{|F_2|\hd |F_k|\}$. So we will consider the
Fubini hypothesis in the form
$$
I_{\cC_{3,x}}=I_{\cC_{2,x}}  \ \ {\rm (Fubini \ hypothesis)}
$$

\smallskip

Now let $X^n\subset \BC\pp{n+2}$ be a submanifold of codimension two
and  satisfy the Fubini hypothesis. What can we say about $X$?

Evident examples for $X$ satisfying the Fubini hypothesis are: the 
intersection of two quadric hypersurfaces, the product of a curve 
with an $(n-1)$-fold having an $n-3$ dimensional family of
lines through a general point  (i.e., a quadric of dimension $n-1$)
or a variety that is a one parameter family of $\pp{n-2}$'s.
Note that to have a meaningful result we should assume $n>2$.

A less evident example is a product of two curves with a $\pp{n-2}$, more
precisely the product of a curve with a variety with a one-dimensional
Gauss image (such varieties are locally the products of curves with
linear spaces). Note that one could not have three curves as
we only have two independent quadrics in the second fundamental form.

We prove

\begin{theorem}[Codimension two Fubini]
\label{codim2fub}
Let $X^n\subset \BC\pp{n+2}$ be an analytic submanifold
with $n>2$. Let $x\in X$ be a general point. If
$$I_{\cC_{2,x}}=I_{\cC_{3,x}} \ \ {\rm (Fubini \ hypothesis)}$$
Then $X$ is one of:
\begin{enumerate}
  \item a complete intersection of two quadric hypersurfaces.
  \item locally the product of a curve with a quadric hypersurface
        $Q^{n-1}\subset\pp{n}$. (I.e., a general point of
$X$ is contained in a $Q^{n-1}\subset X$)
  \item A cone over $Seg(\pp 1\times \pp 2)\subset \pp 5$.
  \item Locally the product of a curve with a variety with a one 
        dimensional Gauss image. In particular, $X$ is locally the product
        of two curves with a $\pp{n-2}$.
  \item Locally the product of a curve with a $\pp{n-1}$,i.e,
a scroll of $\pp{n-1}$'s.
  \item A quadric hypersurface in $\pp{n+1}$.
  \item A linear $\pp n$
\end{enumerate}
\end{theorem}

\bigskip

Under the hypotheses of the theorem $\cC_x$ is the intersection of
(at most) two quadric hypersurfaces. 

The dual 
variety of $X$ is degenerate if and only if none of the quadrics in the ideal of $\cC_x$ are 
smooth. This occurs   in cases 3-7.

\medskip

We expect that our results are valid over $\BR$ in the sense that if one
assumes the same normalizations, the same results hold. However over
$\BR$, there are more cases (e.g., due to the signature of a quadratic
form), although each individual case should be solvable by the methods
of this paper.

\medskip

The meaning of general point here can be made more precise: we assume that 
$|II|_x$ and $|II+F_3|_x$ have base loci having the same number of components 
and dimension of singular sets as all points in some open neighborhood of $x$.
 

\subsection{Related work and problems}

\subsubsection{Rogora's theorem}
By \cite{Lnfac}, Theorem 2, the Fubini hypothesis implies $\cC_{2,x}=\cC_x$ 
(although we do not use this result in our arguments).  Thus a 
generalization  of the problem would be to classify the codimension 
two submanifolds containing (at least) an $(n-3)$-dimensional family 
of lines passing through a general point, or equivalently, the 
codimension two linearly nondegenerate varieties whose Fano variety of 
lines $\BF(X)=\{l\in \BG(\pp 1, \BP V)\mid l\subset X\}$ has dimension 
$2n-4$.  This is a generalization because $\cC_{2,x}$ may have several 
components of dimension $(n-3)$ and the Fubini problem only addresses 
the case when all components are also in $\cC_x$.  
  Now $\tmax\tdim \BF (X)=2n-2$, with equality if and only if 
$X=\pp n$.  The classical Fubini Theorem classifies
the varieties with $\tdim\BF (X)=2n-3$, namely quadric hypersurfaces
and curves of $\pp{n-1}$'s. The next case, where
$\tdim\BF (X)=2n-4$, was solved when $\tcodim (X)>2$ by Rogora
\cite{Rogora}. The only possibilities are one parameter families of quadrics,
two parameter families of $\pp{n-2}$'s or linear sections of
$G(2,5)$, the Grassmannian of $2$-planes in $\BC^5$. The
codimension two case
is partially addressed in this paper:

\begin{corollary} Let $X^n\subset \BP V$ be a projective variety
such that $\tdim \BF_1(X)=2n-4$
and $\cC_{2,x}$ has one component (or such that
$\cC_{2,x}=\cC_{3,x}$). Then unless $X$ is a hypersurface,it is
one of the varieties 1,2,3,4 in the conclusion of theorem \ref{codim2fub}.
\end{corollary}

It is interesting to consider the near counter-example of
a linear projection of $G(2,5)\subset \pp 9$ to a $\pp 8$. In this
case $\cC_{2,x}$ is the union of $\pp 1\times \pp 2$ and a $\pp 3$, but
$\cC_{3,x}=\cC_x=\pp 1\times \pp 2$. This reflects the general principle
that under linear projection from a point, $|II|_x$ loses a quadric
but that quadric shows up multiplied by linear forms in $F_3$.

Similarly, for a  two parameter family of $\pp{n-2}$'s, $\cC_{2,x}$ always has
multiple components.

The case of $\cC_{2,x}$ having  multiple components would be in principle
treatable by the methods of this paper, but one would have to do
a separate calculation for each individual case.
One could study the hypersurface case using the methods of this paper
but it appears one would have to take at least twelve derivatives using the
moving frame to get an answer.

\smallskip

\subsubsection{The Debarre--de Jong conjecture}
Both Debarre and de Jong have conjectured that a smooth hypersurface   
$Z^{m-1}\subset\pp m$ of 
degree $d$ has $\tdim \BF (Z)=2m-3-d$ for $m\geq d$ (the expected
dimension). They observed that by taking
linear sections, it would be sufficient to prove the
conjecture for $d=m$, and moreover 
proved that any potential counter-example $Z$
with a larger space of lines  would have to contain a
hypersurface $X$ (a variety of codimension two in $\pp m$),
with the property that $\BF_1(X)=2m-3-d$, see \cite{beheshti}.

 As an application of our theorem, in section \S\ref{djsect}, we give 
a new proof of this conjecture when $m=6$ (the largest $m$ for which
the conjecture is known to be true), a problem originally
solved by R. Beheshti \cite{beheshti}.

\medskip

\subsubsection{Other generalizations of Fubini}
In codimension one, the Fubini hypothesis implies
that there exists a choice of $F_3$ such that $F_3=0$.
In \cite{Lrigid}, it was shown that a $n$-fold in
$\pp {\binom{n+1}2-1}$ having the expected second fundamental form
(i.e. $|F_2|=S^2T^*_xX$) and admitting a choice of
$F_3$ that is identically zero, must be the quadratic
Veronese embedding of projective space. 
For minimally embedded 
compact Hermitian symmetric spaces (CHSS), something
much stronger is true:
in \cite{Lrigid,Lhss} Fubini's theorem was generalized to all
rank two CHSS in the stronger form that if $\cC_{2,x}$ is the same as that of
a rank two CHSS, then $X$ must be (an open subset of) the corresponding
CHSS. It was then generalized further in \cite{HY} to arbitrary CHSS, requiring
that the Base loci of the fundamental forms coincide. (Roughly speaking,
the $k$-th fundamental form is a  component  of $F_k$ that is
well defined independent of adapted coordinates.)
 
\medskip
 
\subsubsection{An analogue for multi-secant lines?}
Tangent lines are limits of secant lines, and
directions in $\cC_{k,x}$ are limits of $k$-secant lines.
Are there natural analogues of these results
related to $k$-secant lines? For example, much easier than
Fubini's theorem is the fact that a variety $X$  having the
property that any trisecant line is contained in $X$ is
either a quadric or a linear space.


\subsection{Outline of the proof}\label{sec:outline}  

We know of two proofs of Fubini's result (Theorem \ref{thm:fubini}). 
One can either reduce to
the surface case by taking a general $\pp 3$-section and then
prove the theorem for surfaces (which follows because a surface
having two distinct lines through a general point is necessarily
a quadric) or by reducing the frame bundle of an unknown variety
satisfying the Fubini hypothesis to the reduced frame bundle of
a quadric hypersurface. Any proof of the codimension two Fubini
theorem must necessarily be more complicated because for quadric
hypersurfaces (the codimension 1 case), there is only a discrete invariant  
(the rank), but for pencils of quadrics (the codimension 2 case) there are 
moduli. Thus a moving frames proof would have to reduce to a Frobenius system 
on the frame
bundle (i.e., one whose solutions were parametrized by a fixed
number of constants).  For a linear section argument, one needs to
be sure that the sections cannot be coming from a more complicated
variety (since the sections will not all be isomorphic).  

Moreover, not only do the expected answers have moduli, the
possible second fundamental forms do as well (as they too are
pencils of quadrics), whereas in the original Frobenius theorem
there was only the discrete invariant of rank.  Our proof combines
methods of both proofs of Fubini's theorem.

If a variety satisfies Fubini's hypothesis,
then so will any general linear section.  For most cases
we  prove Theorem 
\ref{codim2fub} for $n=3$  and then use the fact that any general $\PP^5$
section of $X^n$ is of the type found in the $n=3$ analysis to characterize
these varieties. In other cases we just argue directly in $n$ dimensions.
For the generic $\cC_{2,x}$ both methods work equally well. Here
are the possible cases:
\begin{itemize}
  \item[(1)] Whenever a general linear section of a variety is a complete 
        intersection cut out by varieties of degrees $d_1\hd d_s$, then the 
        original variety must also be a complete intersection 
        cut out by varieties of degrees $d_1\hd d_s$.  
  \item[(2)]   Here we  
        prove the result directly for arbitrary $n$.
  \item[(3)] The only variety whose general $\pp 5$ section  is 
        $\pp 1\times \pp 2$ is a cone over $\pp 1\times \pp 2$.   
  \item[(4)] If a general linear section of $X$ is locally the product of a 
        curve with a variety with a one dimensional Gauss image, then $X$ 
        will have that property as well.  
  \item[(5-7)] These cases are degenerate and covered by the original
        (codimension one) Fubini theorem remarks.
\end{itemize}

Now consider the case $n=3$. We have the following possibilities for
$|II|_{X,x}$ where $x\in X$ is a general point. If it consists of a single 
quadric, either the quadric has rank greater than one and $X$ is a quadric 
hypersurface in a $\PP^{n-1}$, or else $X$ is a curve of $\pp{2}$'s. If it is 
a pencil, then, as 
explained in \cite{HP}, there are seven possibilities for the pencil, as 
characterized by the base loci:
If the pencil contains a smooth conic, then
the base locus $\cC_{2,x}$ consists of four points
(counted with multiplicity)
in $\pp 2=\BP T_xX$. The cases are:
  (i) four distinct points;
 (ii) two double points;
(iii) a double point and two distinct points;
 (iv) a single four-fold point;
  (v) a triple point and a distinct point.
The sixth and seventh cases arise when the pencil contains no smooth 
quadrics.  Equivalently, the dual is degenerate.
The seven cases are analyzed in Subsections \S\ref{sec:Case1}--\ref{dualandg}.

\medskip

 \noindent{\it Acknowledgments.} We would like to thank I. Coskun and F. Zak for useful remarks, and in particular we thank Coskun for providing
us with theorem 3.1.

  
\section{Moving frames}
\label{sec:frames}
 
We use  notation for the moving frame and
differential invariants as in
\cite{IvL}.
We use index ranges
\begin{align*}
&1\leq a,b,..,e\leq n\\
&n+1\leq u,v\leq n+2\\
&0\leq A,B\leq n+2.
\end{align*}

\noindent{\bf NOTE}: In calculations we will use the convention that indices
$a,b$ are {\it not} to be summed over unless explicitly specified
but use the summation convention for all other indices.

We work on the open subset of a codimension 2 submanifold 
$X^n \subset \CP^{n+2}$ consisting of general points and
slightly abuse notation by calling it $X$.

The bundle of first order adapted frames $\cF^1_X$ for a 
submanifold $X^n\subset \pp{n+2}=\BP V$ is the set  of ordered bases
$g=(e_0\hd e_{n+1})$ of $V$ such that $[e_0]\in X$ and
the affine tangent space
$\hat T_xX$ is the span of $e_0\hd e_n$. It is a bundle over 
$X$ and the Maurer-Cartan form  $\w=(\o^A_B)=g\inv dg$ 
of $GL(V)$ pulls back to give
forms on $\cF^1_X$.  
We write $g=(g^A_B)\in GL(V)$.

The first order adaption forces 
\begin{displaymath}
   \w^u_0  = 0 \, .
\end{displaymath}
Differentiating these equations produces
\begin{equation}
\label{eqn:q}
  \w^u_a = q^u_{ae} \w^e_0 \, ,
\end{equation}
for symmetric functions $q^u_{ab} = q^u_{ba}$.
A moving frame  definition of 
  the second fundamental form
$F_2=II_{X}\in\Gamma(X,S^2T^*X\ot NX)$
is obtained by pushing down $\ooo ue\ot \oo e\ot e_u
\in\Gamma(\cF^1_X,\pi^*(S^2T^*X\ot NX))$ down to $X$.  We denote the Fubini cubic 
by $F_3=r^u_{efg}\oo e\oo f\oo g\ot e_u\in\Gamma(\cF^1_X,\pi^*(S^3T^*X\ot NX))$
where the coefficients
  $r^u_{abc}$ of $F_3$ are defined by 
\begin{equation}
\label{eqn:F3_def}
  r^u_{abc} \w^c_0 = -d q^u_{ab} - q^u_{ab} \w^0_0 - q^v_{ab} \w^u_v 
                   + q^u_{ae} \w^e_b + q^u_{be} \w^e_a \, .
\end{equation}
See \cite{IvL}, Chapter 3 for details.

We now add the Fubini  hypothesis  that
$|F_3|\subset |II\circ T^*|$ on the coefficients of $F_3$:
\begin{displaymath}
  r^u_{abc} = \gS_{abc} \, \rho^u_{av} q^v_{bc} \, .
\end{displaymath}
The notation $\gS$ denotes cyclic summation on the indices.

The two degenerate cases (vi,vii) for $n=3$ mentioned in Subsection 
\ref{sec:outline} have the normal forms $\{\oo 1\oo 2, \oo 1\oo 3\}$ and 
$\{(\oo 1)^2,(\oo 2)^2\}$. These cases are respectively treated in 
subsections \S\ref{dualnotg} and \S\ref{dualandg}.


\subsection{Case (i):  {\boldmath{$n=3$, and $\cC_{2,x}$}} is linearly 
                    nondegenerate and smooth}
\label{sec:Case1}
Here we begin our seven part analysis of the case that $n=3$ and $|II|_{X,x}$ 
contains a pencil of quadrics.  When the base locus $\cC_{2,x}$ contains four 
distinct points we may normalize the $q^u_{ab}$ so that 
\begin{equation}
\label{eqn:Case1_q}
  q^{n+1}_{ab} = \d_{ab} 
  \qquad\hbox{ and }\qquad 
  q^{n+2}_{ab} = \l_a \d_{ab} \,  
\end{equation}
for pairwise distinct functions $\l_a$.  
We will see that $X$ is the intersection of two quadrics.  To do so it is 
sufficient, by Theorem 4.28 of \cite{Lci}, to show that the coefficients of 
$F_4$ and $F_5$ satisfy
\begin{eqnarray}
\label{eqn:F4}
  r^u_{abcd}& = & \mathfrak{S}_{abc} \, \sigma^u_{vw} \, 
                                       q^v_{ab} \, q^w_{cd} \ + \
                  \mathfrak{S}_{abcd} \, \rho^u_{av} \, r^v_{bcd} \, , \\
\label{eqn:F5}
  r^u_{abcde} & = & \gS_{abcde} \, \left( 
                   \rho^u_{av} \, r^v_{bcde} + 
                   \sigma^u_{vw} \, ( q^v_{ab} \, r^w_{cde}
                                    + q^v_{ac} \, r^w_{ebd} )
                   \right) \, . 
\end{eqnarray}
(Although we will not need to use this in our calculations, this
will serve as a useful guide.)
Here $\sigma^u_{vw} = \sigma^u_{wv}$.  Recall that these coefficients are 
defined by 
\begin{eqnarray}
\label{eqn:F4_def}
  r^u_{abcd} \w^d_0 & = & - d r^u_{abc} - 2 r^u_{abc} \w^0_0 
                        - r^v_{abc} \w^u_v \\
  \nonumber       &   & + \gS_{abc} \big( 
                                            r^u_{abe} \w^e_c + q^u_{ab} w^0_c 
                                          - q^u_{ae} q^v_{bc} \w^e_v
                                    \big) \, , \\
\label{eqn:F5_def}
  r^u_{abcde} \w^e_0 & = & - d r^u_{abcd} - 3 r^u_{abcd} \w^0_0 
                         - r^v_{abcd} \w^u_v \\
  \nonumber        &   & - \gS_{abc} 
                           \big\{
                             ( r^u_{abe} q^v_{cd} + r^u_{ade} q^v_{bc} ) \w^e_v
                           + ( q^u_{ab} q^v_{cd} + q^u_{ad} q^v_{bc} ) \w^0_v
                           \big\}  \\
  \nonumber        &   & + \gS_{abde}
                           \big\{
                             r^u_{abce} \w^e_d + 2 r^u_{abc} \w^0_d
                             - q^u_{ae} r^v_{bcd} \w^e_v 
                           \big\} \, .
\end{eqnarray}

We will use the notation 
\begin{displaymath}
  \l_{ab} := \l_a - \l_b \not= 0 \, .
\end{displaymath} 
Note that (\ref{eqn:q}) gives us
\begin{equation}
\label{eqn:2ndOrder}
  \w^{n+1}_a = \w^a_0
  \qquad\hbox{ and }\qquad
  \w^{n+2}_a = \l_a \w^a_0 \, .
\end{equation}
Recall our convention that there is no sum on $a$ in the last equation.   

Assume Fubini's hypothesis holds.  For a suitable choice of $g^0_a , g^a_v$, 
the transformation $e_u \mapsto e_u + g^a_u e_a$ and 
$e_a \mapsto e_a + g^0_a e_0$ further refines the frames so that 
$\rho^4_{av} = 0 = \rho^5_{a5}$ and $\rho^5_{a4} = \rho_a$.  
Now (\ref{eqn:F3_def}) implies Fubini's hypothesis holds 
on our reduced frame bundle if and only if
\begin{eqnarray}
  \label{eqn:fubini1}
  0 & = & \w^a_b + \w^b_a  \hbox{\hspace{2in}} (a \not= b) \\
  \label{eqn:fubini2}
  0 & = & - \w^0_0 - w^{n+1}_{n+1} - \l_a \w^{n+1}_{n+2} + 2 \w^a_a \\
  \label{eqn:fubini3}
  \r_a \w^b_0 + \r_b \w^a_0 & = & \l_a \w^a_b + \l_b \w^b_a 
    \hbox{\hspace{1.7in}} (a \not= b) \\ 
  \label{eqn:fubini4}
  2 \r_a \w^a_0 + \r_e \w^e_0 & = & - d\l_a - \l_a \w^0_0 - \w^{n+2}_{n+1}
                               - \l_a \w^{n+2}_{n+2} + 2 \l_a \w^a_a \, .
\end{eqnarray}
(The first two equations come from $u = n+1$, and the last two from $u=n+2$.)
\subsubsection{Determination of $F_4$}
\label{ssec:stage1}
Differentiating (\ref{eqn:fubini1}) produces functions $C^a_e$ and $E^a_e$ so 
that
\begin{eqnarray*}
    \w^a_5 & = & C^a_e \w^e_0 \\
    \w^0_a - \w^a_4 - \r_a \w^4_5 & = & E^a_e \w^e_0 \, .
\end{eqnarray*}
The functions $E^a_e$ satisfy the relations 
\begin{displaymath}
\renewcommand{\arraystretch}{1.3}
  (E^a_e) = \left( \begin{array}{ccc}
              E^1_1 & \l_3 C^1_2 & \l_2 C^1_3 \\
              \l_3 C^2_1 & E^2_2 & \l_1 C^2_3 \\
              \l_2 C^3_1 & \l_1 C^3_2 & E^3_3 
            \end{array}\right) \, ,
\end{displaymath}
and 
\begin{displaymath}
  C^a_a \l_b - C^b_b \l_a = E^a_a - E^b_b \, .
\end{displaymath}
The last is a set of $\binom n2 = 3$ linear equations for the $2n = 6$ 
unknowns $C^a_a,E^b_b$ but the system has rank $2n-4 = 2$ so there is a 
$4$-dimensional space of solutions. We may parameterize the solutions as 
follows by introducing new variables $R,S,T$:
\begin{equation}
\label{eqn:CE}
   C^a_a = R \l_a + S \qquad \hbox{ and } \qquad
   E^a_a = -S \l_a + T \, .
\end{equation}
The derivative of (\ref{eqn:fubini2}) forces the off-diagonal terms of $C$, 
and therefore $E$ as well, to vanish.  Whence
\begin{eqnarray}
  \label{eqn:wa2_i}
    \w^a_5 & = & C^a_a \w^a_0 \\
  \label{eqn:w0a_i}
    \w^0_a - \w^a_4 - \r_a \w^4_5 & = & E^a_a \w^a_0 \, .
\end{eqnarray}
We may use $g^0_4,g^0_5$ to normalize $S,T=0 \ \Longrightarrow \ E^a_a = 0$.

Making use of the identities derived thus far, a computation of 
(\ref{eqn:F4_def}) in the $u=n+1=4$ case yields 
\begin{displaymath}
  r^4_{abcd} = \gS_{abc} \sigma^4_{uv} q^u_{ab} q^v_{cd} \, ,
\end{displaymath}
with $\sigma^4_{4,4} = \sigma^4_{4,5} =   \sigma^4_{5,4}=0$ and 
$\sigma^4_{5,5} = -R$.  Taking into account the normalizations of 
$\rho$, this gives us the $u=n+1=4$ half of (\ref{eqn:F4}). 

Next we differentiate (\ref{eqn:fubini3}) and obtain functions $F^a_e, G^a_e$
such that
\begin{eqnarray*}
  w^a_4 + 2 \r_a w^4_5 & = & F^a_e \w^e_0 \\
  d \r_a - \r_a ( 2 w^4_4 + 3 \l_a \w^4_5 - w^5_5 )
         + \sum_e \r_e \w^a_e & = & G^a_e \w^e_0 \, .
\end{eqnarray*}
Additionally, the functions $G^a_e$ satisfy 
\begin{displaymath}
\renewcommand{\arraystretch}{1.3}
   (G^a_e) = \left( \begin{array}{ccc}
               G^1_1 & \l_{31} F^1_2 & \l_{21} F^1_3 \\
               \l_{32} F^2_1 & G^2_2 & \l_{12} F^2_3 \\
               \l_{23} F^3_1 & \l_{13} F^3_2 & G^3_3 
             \end{array} \right) \, ,
\end{displaymath}
and 
\begin{displaymath}
  ( F^a_a - \l_a{}^2 R ) \l_b - ( F^b_b - \l_b{}^2 R ) \l_a 
   = 
  \Big(
    G^a_a + \l_a F^a_a
  \Big)
  - 
  \Big(
    G^b_b + \l_b F^b_b
  \Big) \, .
\end{displaymath}
As above for (\ref{eqn:CE}), this is a corank three system and, introducing 
new variables $U,V,W$ gives,
\begin{equation}
\label{eqn:G}
\begin{array}{rcl}
   F^a_a - \l_a{}^2 R & = & U \l_a + V \\
   \displaystyle
   G^a_a + \l_a F^a_a & = & -V \l_a + W \, .
\end{array}
\end{equation}
The derivative of (\ref{eqn:fubini4}) forces the off-diagonal entries of 
$F$ (and therefore $G$, as well) to vanish.  With an application of 
(\ref{eqn:w0a_i}) we have  
\begin{eqnarray}
\label{eqn:F_i}
  \w^0_a + \rho_a \w^4_5 & = & F^a_a \w^a_0 \, , \\
\label{eqn:drho_dst}
  d \rho_a 
  - \r_a ( 2 w^4_4 + 3 \l_a \w^4_5 - w^5_5 )
  + \sum_e \r_e \w^a_e
  & = & G^a_a w^a_0 \, .
\end{eqnarray}

Now a computation of (\ref{eqn:F4_def}) in the $u=n+2=5$ case yields 
 \begin{displaymath}
  r^5_{abcd} = \gS_{abc} \sigma^5_{uv} q^u_{ab} q^v_{cd} \, ,
\end{displaymath}
with 
$\sigma^5_{4,4} = -W$, 
$\sigma^5_{4,5} = V = \sigma^5_{5,4}$ and 
$\sigma^5_{5,5} = U $.  In particular, (\ref{eqn:F4}) holds.

\subsubsection{Determination of $F_5$}
\label{ssec:stage2}
It remains to verify (\ref{eqn:F5}).  These coefficients are given by 
(\ref{eqn:F5_def}) which requires that we compute $-d r^u_{abcd}$.  In 
particular, we   need expressions for $dR$,   $dU$, $dV$ and $dW$.
We   obtain information on the first three differentials by differentiating 
the expressions 
\begin{eqnarray}
\label{eqn:wa2_ii}
  \w^a_5 -   R \l_a   \w^a_0 & = & 0 \, \quad \hbox{ and } \\
\label{eqn:w0a_ii}
  \w^0_a - \w^a_4 - \rho_a w^4_5    & = & 0 \, ,
\end{eqnarray}
which are consequences of (\ref{eqn:CE},\ref{eqn:wa2_i},\ref{eqn:w0a_i}).
In particular, we find
\begin{equation}\label{eqn:dsigma1}\renewcommand{\arraystretch}{1.3}
\begin{array}{rcl}
  dR & = & R ( 2 \w^5_5 - \w^0_0 - \w^4_4 ) 
           +  U  \w^4_5 \\
  0 & = & \w^0_5 +    
           +  V  \w^4_5 + R \w^5_4 \\
  0 & = & -2 \w^0_4 +  W \w^4_5  \, .
\end{array}\end{equation}
The first and second expressions   are derived from the derivative of 
(\ref{eqn:wa2_ii}), and the third from (\ref{eqn:w0a_ii}).

Next, (\ref{eqn:G},\ref{eqn:F_i},\ref{eqn:drho_dst}) give us
\begin{eqnarray}
\label{eqn:F_ii}
   0 & = & 
   \w^0_a + \rho_a \w^4_4 - ( R \l_a{}^2 + U \l_a + V ) \w^a_0 \\
\label{eqn:G_ii}
   0 & = & 
   \left(
     R \l_a{}^3 +  U  \l_a{}^2 +   2 V   \l_a - W
   \right) \w^a_0 \\
   \nonumber & &
   + \, d \rho_a - \rho_a ( 2 w^4_4 + 3 \l_a \w^4_5 - \w^5_5 )
            + \sum_e \r_e \w^a_e \, .
\end{eqnarray}
Differentiating (\ref{eqn:F_ii}) provides expressions for $dU$ and $dV$; 
$dW$ is given by (\ref{eqn:G_ii}).
Summing over the index $e \in \{ 1, \ldots , n \}$, we have
\begin{equation}\label{eqn:dsigma2}\renewcommand{\arraystretch}{1.3}
\begin{array}{rcl}
  dU & = & -\w^0_5 + R \rho_e \w^e_0 + U ( \w^5_5 - \w^0_0 )
           +   3V  \w^4_5 + 2R \w^5_4 \\
  dV & = & -\w^0_4 - R \l_e \rho_e \w^e_0 + V( w^4_4 - \w^0_0 )
           - W \w^4_5 + U \w^5_4 \\
  dW & = & ( 2 R \l_e{}^2 + 2  U  \l_e + 2 V - T ) \rho_e \w^e_0 
           - 4 \rho_e{}^2 \w^4_5 \\
     &   & + W ( 2 \w^4_4 - \w^0_0 - \w^5_5 )
           -  2V  \w^5_4 \, .
\end{array}\end{equation}
Now a computation of (\ref{eqn:F5_def}) reveals that the coefficients of 
$F_5$ are indeed of the form (\ref{eqn:F5}), and $X$ must be a complete 
intersection in the case of distinct eigenvalues for the second quadric in 
$II$.

\medskip Note that one can avoid the use of \cite{Lci}, Theorem 4.28 as 
follows:  Differentiating (\ref{eqn:dsigma1},\ref{eqn:dsigma2}) yields no 
additional relations and we may make the following observation.  Let $\C^3_\l$ 
and $\C^3_\rho$ denote two copies of $\C^3$ with coordinates 
$\l = (\l_1,\l_2,\l_3)$ and $\rho = (\rho_1,\rho_2,\rho_3)$, respectively.
Denote the coordinates of $\C^4_\sigma$ by $(R,U,V,W)$.  Let 
$M = \{ (\l_1,\l_2,\l_3) \in \C^3_\l \ : 
                \ \l_a \not= \l_b \hbox{ whenever } a \not=b \}$, and 
$\Sigma = GL_{n+3}\C \times M \times \C^3_\rho \times \C^4_\sigma$
(here $n=3$).  Then 
the system given by the equations $\{ \w^u_0 = 0 \}$ and 
(\ref{eqn:2ndOrder},\ref{eqn:fubini1},\ref{eqn:fubini2},\ref{eqn:fubini3},\ref{eqn:fubini4},\ref{eqn:wa2_ii},\ref{eqn:w0a_ii},\ref{eqn:dsigma1},\ref{eqn:F_ii},\ref{eqn:G_ii},\ref{eqn:dsigma2}) is Frobenius.  Note that 
$\hbox{dim}_\C \Sigma = 46$, and that the system consists of 36 independent 
equations.  So the maximal integral submanifolds are of dimension 10 and may 
be identified with the graphs of those the natural maps 
$\mathcal{F} \to M \times \C^3_\rho \times \C^4_\sigma$, where 
$\mathcal{F} \subset GL_{n+3}\C$ is a sub-bundle of the adapted frame bundle
over a smooth variety of codimension 2 which satisfies Fubini's hypothesis
(and with distinct eigenvalues $\l_a$).
In fact the resulting integral manifolds have ideal generated by
\begin{align*}
&\xx 0\xx{4} - \sum_a(\xx a)^2  
+R(\xx{5})^2,\\
&\xx 0\xx{5} - \sum_a\l_a(\xx a)^2
 - \sum_a\rho_a\xx a\xx{4} +W(\xx{4})^2 -V\xx{4}\xx{5}
 -U (\xx{5})^2.
\end{align*}

{\it Remark.} This computation is easily generalized to arbitrary $n$.
In particular, suppose the second quadric in $II$ may be normalized as in 
(\ref{eqn:Case1_q}).  Additionally assume that there exists as least two
distinct eigenvalues $\l_a$, and that no eigenvalue occurs with multiplicity
$n-1$.  (In the case $n=3$ this is equivalent to hypothesis of three distinct
eigenvalues.)  The analogous calculation for $n>3$ shows that each numbered 
equation in this section holds when the indices $(4,5)$ are replaced with 
$(n+1,n+2)$.  Again we have a Frobenius
system whose integral manifolds have ideal generated by
\begin{align*}
&\xx 0\xx{n+1} - \sum_a(\xx a)^2  
+R(\xx{n+2})^2,\\
&\xx 0\xx{n+2} - \sum_a\l_a(\xx a)^2
 - \sum_a\rho_a\xx a\xx{n+1} +W(\xx{n+1})^2 -V\xx{n+1}\xx{n+2}
 -U (\xx{n+2})^2 \, .
\end{align*}


\subsection{Case (ii).}
\label{sec:Case2}
Fix $n\geq 3$ and assume we are in Case (ii).  As in Case (i) we assume the $q^u_{ab}$ may be normalized so that 
\begin{displaymath}
  q^{n+1}_{ab} = \d_{ab} 
  \qquad\hbox{ and }\qquad 
  q^{n+2}_{ab} = \l_a \d_{ab} \, .
\end{displaymath}
Additionally, assume a $1,n-1$ split of the eigenvalues: 
$\l_1 \not= \l_2 = \cdots = \l_n$.  (When $n=3$, this is the case that 
$\cC_{2,x}$ contains two points, each counted with multiplicity 2.)

A second normalization puts the coefficients in the form $\l_1=0$ and 
$\l_{\a}=\l \not=0$ for
$2\leq \a \leq n$. Then \eqref{eqn:fubini1},\eqref{eqn:fubini3} imply
$$
\w^\a_1 = - \ooo 1\a = \frac{\r_1}{\l} \oo\a
$$
which in turn implies the hyperplane distribution $\{\oo 1\}\upperp$ 
is integrable.
Since any line field is integrable as well, we see that $X$ is
locally a product $C\times Y$. But now $II_Y$ consists of a single
quadric of rank greater than one, so by
e.g. \cite{IvL}, Cor. 3.5.7, $Y$ is a hypersurface in some $\pp{n-1}$.
Fubini's hypothesis also holds for $Y$ so that it must be a quadric 
hypersurface by Theorem \ref{thm:fubini}.  This places us in Case (2) of 
Theorem \ref{codim2fub}.


\subsection{Case (iii).} 
\label{sec:Case3}  This is the case that $\cC_{2,x}$ consists of 3 points, 
one with multiplicity 2.  We may normalize the second fundamental form as 
follows 
\begin{displaymath}
  \left( q^4_{ab} \right) = 
  \left(\begin{array}{ccc}
    0 & 0 & 0 \\
    0 & 1 & 0 \\
    0 & 0 & \l
  \end{array}\right) 
  \qquad\hbox{ and }\qquad 
  \left( q^5_{ab} \right) = 
  \left( \begin{array}{ccc}
    0 & 1 & 0 \\
    1 & 0 & 0 \\
    0 & 0 & 1
  \end{array}\right)  
\end{displaymath}
for some function $\l \not=0$.  In particular, (\ref{eqn:q}) gives us
\begin{equation}
\label{eqn:I_2ndOrder}
  \renewcommand{\arraystretch}{1.3}
  \begin{array}{rclcrclcrcl}
    \w^4_1 & = & 0 & & \w^4_2 & = & \w^2_0 & & \w^4_3 & = & \l \w^3_0 \\
    \w^5_1 & = & \w^2_0 & & \w^5_2 & = & \w^1_0 & & \w^5_3 & = & \w^3_0 \, .
  \end{array}
\end{equation}

Assume Fubini's hypothesis holds.  As in \S\ref{sec:Case1} a suitable choice
of $g^0_a , g^a_v$ allows us to normalize $\rho$.  In this case we may refine
the framing so that $\r^5_{av} = 0 = \r^4_{a4}$ and $\r^4_{a5} = \r_a$.
(Contrast with \S\ref{sec:Case1},\ref{sec:Case2} where $\r_a = \r^5_{a4}$.)
Computations with (\ref{eqn:F3_def}) produce
\begin{eqnarray*}
  (u,a,b) = (5,1,1) & \Longrightarrow & \w^2_1 = 0 \\
  (u,a,b) = (4,3,1) & \Longrightarrow & \l \w^3_1 = \r_3 \w^2 \\ 
  (u,a,b) = (5,3,1) & \Longrightarrow & \w^2_3 + \w^3_1 = 0 \, .
\end{eqnarray*}
The last two equations tell us that $\w^2_3 \equiv 0$ mod $\w^2_0$.  Along 
with the first equation above, this implies the hyperplane distribution 
$\{ \w^2_0 \}^\perp$ is integrable.  As in \S\ref{sec:Case2} $X$ is locally 
the product of a curve and surface $Y$.  In this case 
$$
 II_Y  = (\ooo 41\oo 1+\ooo 43\oo 3)\ot e_4+
 (\ooo 51\oo 1+\ooo 53\oo 3)\ot e_5+
(\ooo 21\oo 1+\ooo 23\oo 3)\ot e_2 \tmod \oo 2 ,
$$ 
so 
$|II_Y|= \{ (\w^3_0)^2 \}$.
Hence the Gauss map of $Y$ is degenerate and we are in Case (4) of Theorem 
\ref{codim2fub}.  (Note that when $n=3$, we can say more as $Y$ is either a 
cone over a curve, or the tangential variety of a curve.  (Cf. \cite{CS10}, 
p.105; or \cite{IvL}, Thm.3.4.6.))


\subsection{Case (iv).}
\label{sec:Case4}
Here $\cC_{2,x}$ contains a single point of multiplicity 4.  The second 
fundamental form may be normalized as in \S\ref{sec:Case3}, but with $\l=0$.
Again,
\begin{eqnarray*}
  (u,a,b) = (4,3,1) & \Longrightarrow & \r_3 = 0 \\
(u,a,b) = (5,1,1) & \Longrightarrow & \w^2_1 = 0 \\
  (u,a,b) = (4,1,2) & \Longrightarrow & -\w^4_5 = 2 \r_2 \w^2_0 \\
  (u,a,b) = (4,3,3) & \Longrightarrow & -\w^4_5 = \r_2 \w^2_0
\end{eqnarray*}
Those last two equations imply $\r_2 = 0$.  Now 
\begin{eqnarray*}
  (u,a,b) = (4,2,3) & \Longrightarrow & \w^2_3 = 0\, ,
\end{eqnarray*}
and $\{ \w^2_0 \}^\perp$ is again integrable
and again $|II_Y|=\{(\oo 3)^2\}$ and we are in Case (4) as above.


\subsection{Case (v).}
Here $\cC_{2,x}$ consists a triple point, and a singleton.  We may normalize 
the second fundamental form as follows 
\begin{displaymath}
  \left( q^4_{ab} \right) = 
  \left(\begin{array}{ccc}
    0 & 0 & \l \\
    0 & \l & 1 \\
    \l & 1 & 0
  \end{array}\right) 
  \qquad\hbox{ and }\qquad 
  \left( q^5_{ab} \right) = 
  \left( \begin{array}{ccc}
    0 & 0 & 1 \\
    0 & 0 & 0 \\
    1 & 0 & 1
  \end{array}\right)  
\end{displaymath}
for some function $\l \not=0$.  The Fubini cubic may be normalized so that 
$\r^4_{av} = 0 = \r^5_{a5}$ and $\r_a := \r^5_{a4}$.  Now computations of 
(\ref{eqn:q}) yield
\begin{eqnarray*}
  (u,a,b) = (4,1,1) & \Longrightarrow & \w^3_1 = 0 \\
  (u,a,b) = (5,1,1) & \Longrightarrow & \w^3_1 = \l \r_1 \w^3_0 \\
  (u,a,b) = (5,1,2) & \Longrightarrow & \w^3_2 = \l \r_1 \w^2_0
                                               + ( \r_1 + \l \r_2 ) \w^3_0 \, .
\end{eqnarray*}
In particular, $\r_1 = 0$, and $\w^3_2 \equiv 0$ mod $\w^3_0$.  It follows that
the hyperplane distribution $\{ \w^3_0 \}^\perp$ is integrable
and $|II_Y|=\{(\oo 2)^2\}$ and again we are in Case (4).


\subsection{Degenerate dual and nondegenerate Gauss map case}
\label{sec:deg_dual}\label{dualnotg}
Here there is a unique pencil of
quadrics up to equivalence satisfying
the hypotheses: we may normalize $|II_X|=\{\oo 1\oo 3,\oo 2\oo 3\}$.  Now
the hypothesis on $F_3$ allows us to reduce the frame bundle on $X$ to a 
sub-bundle upon which the Maurer-Cartan forms pull-back to satisfy the same 
relations as those satisfied by the Maurer-Cartan forms on the frame bundle of 
$\pp 1 \times \pp 2$ in its Segre embedding.  In particular, both bundles are 
integral manifolds of a Frobenius system defined by left-invariant 1-forms on 
$GL_6\C$.  Hence, $X$ is (projectively equivalent to an open subset of) 
$\hbox{Seg}(\pp 1\times \pp 2)$.


\subsection{Degenerate dual and rank two Gauss map case}\label{dualandg}
Here we may normalize $|II_X|=\{(\oo 1)^2, (\oo 2)^2\}$. This case
also reduces to Case (4) and the calculation is even easier than
the above cases.


\section{Proof of the Debarre-de Jong  conjecture for degree six hypersurfaces}
\label{djsect}
We need to show no smooth hypersurface $Z^5 \subset \PP^6$ of degree six can
contain a codimension two subvariety with $4$-dimensional
Fano variety of lines. Our proof will use general results
of \cite{beheshti} but avoid the case by case
study in section 4.2 of \cite{beheshti}. Our proof may be useful
in either proving the degree seven case or as a guide to
potential counter-examples for all higher degrees.

We have classified $4$-folds in $\pp 6$ with
a $4$ dimensional Fano variety of lines when 
$\cC_{2,x}$ has one component. When it has several components,
they must be curves of degrees one, two or three (and non-planar
in the last case; see \cite{HP}, p.307, case (x)), and therefore
rational, so by
\cite{beheshti}, Theorem 2.1, cannot lie in a smooth $5$-fold hypersurface.

The intersection of two quadrics is ruled out by degree considerations.
A variety that is locally the product of a curve with a $3$-quadric
is ruled out again because $\cC_x$ contains a plane conic.
And a cone over $\pp 1\times \pp 2$ is similarly ruled out.

It remains to deal with the case when $X^4$ is contained in a $\pp 5$,
and therefore would be a hyperplane section of some
smooth counter-example $Z^5\subset \pp 6$.
In fact we know it would have to be a singular hyperplane section, i.e.,
$X=Z\cap \tilde T_zZ$ for some $z\in Z$.
(Note that $X$ must be uniruled by lines, as if the lines of
$X$ passed through a proper subvariety, a component of that
subvariety would have to be a $\pp 3$ which cannot be in
a smooth $5$-fold.) Thanks to Zak's theorem on
tangencies \cite{zak}
 we know such $X$ has at most isolated singularities (the
singular points are the other points of $Z$ tangent to
the hyperplane $\tilde T_zZ$).
Note also that $X$ cannot be a cone as then it would support at most
a $3$ dimensional family of lines (unless it were a $\pp 4$).
To finish we appeal to a result supplied
to us by I.  Coskun (personal communication)
which he believes to be \lq\lq known to the experts\rq\rq :

\begin{theorem} Let $X^n\subset\pp{n+1}$   be a hypersurface of 
degree at least $n+1$. Suppose $X$ has only isolated singularities and $X$ 
is not a cone. Then $X$ is not covered by lines. 
\end{theorem}

\begin{proof}
Since $X$ is not a cone, a general line on $X$ has a well defined
normal bundle $N_{L/X}$ over $L$.  This bundle if of rank $n-1$ 
and fits into the exact sequence
$$
0\ra N_{L/X} \ra N_{L/\pp{n+1}}\ra N_{X/\pp{n+1}}\mid_L\ra 0 \, 
$$
Note that the second term is just $\cO(1)^{\op{n}}$ and the last is $\cO(d)$ 
where $d\geq n+1$ is the degree of $X$.  Write 
$N_{L/X}=\cO(a_1)\op \cdots \op \cO(a_{n-1})$ (by the Segre-Grothendieck 
splitting theorem), with $a_1\leq a_2\leq \cdots \leq a_{n-1}$.
Then we see $a_1<0$ which means the deformations of $L$ cannot
cover $X$, a contradiction.
\end{proof}


\end{document}